\documentclass[oneside,12pt]{amsart}
\usepackage{amssymb, amscd}
\usepackage[all]{xy}

\theoremstyle{plain}
\newtheorem{theorem}{Theorem}[section]
\newtheorem{lemma}[theorem]{Lemma}

\newtheorem*{theorem*}{}
\newtheorem{proposition}[theorem]{Proposition}
\newtheorem{corollary}[theorem]{Corollary}

\theoremstyle{definition}

\newtheorem{question}[theorem]{Question}
\newtheorem{definition}[theorem]{Definition}

\theoremstyle{remark}

\DeclareMathOperator{\Z}{\mathbb Z}
\DeclareMathOperator{\Q}{\mathbb Q}

\begin{document}

\title{Almost Cyclic Groups}
\author{Bruce Ikenaga\\
   Department of Mathematics\\
   Millersville University\\
   Millersville, PA 17551}
\maketitle

\section{Introduction.}

\begin{definition} A group $G$ is {\bf almost cyclic} if
there exists an element $x \in G$ such that every $y \in G$ is conjugate
to a power of $x$. Thus, for all $y \in G$, there is an element $a \in
G$ such that
$$aya^{-1} = x^n \quad\hbox{for some}\quad n \in \Z.$$

$x$ is called a {\bf conjugate generator} for $G$.
\end{definition}

If a conjugate generator has finite order $n$, then $n$ is an exponent
for the group. On the other hand, if a conjugate generator has infinite
order, then the group is torsion-free. Thus, the class of almost cyclic
groups divides naturally into the torsion almost cyclic groups and the
torsion-free almost cyclic groups.

If $G$ is a finitely generated torsion almost cyclic group, then $G$ is
a counterexample to the Burnside conjecture (of which many have been
known for some time). However, it suggests that constructing examples of
such groups will be difficult. For certain exponents, known results
apply. For example, suppose $G$ is almost cyclic with conjugate
generator $x$ of order $2$. Then every nontrivial element of $G$ has
order $2$, so $G$ is abelian. This implies (as will be shown below) that
$G$ is isomorphic to $\Z_2$.

Every cyclic group is almost cyclic. Fortunately, there are almost
cyclic groups which are not cyclic The key problem is to determine when
an almost cyclic group is cyclic; consequently, many of the questions
considered in this paper have the following form. Let {\bf P} be a
group-theoretic property. Is an almost cyclic group having {\bf P}
necessarily cyclic?

Almost cyclic groups were discussed by W. Ziller in an unpublished paper
\cite{ziller} in connection with the existence of closed geodesics on
manifolds. He conjectures (in the terminology of the present paper) that
every finitely presented almost cyclic group is cyclic. He notes further
that if this is true, then on any compact manifold $M$ with $\pi_1(M)
\ne \{1\}, \Z, \Z_n$, one obtains at least two closed geodesics from the
fundamental group. Ziller's conjecture appears to be difficult to
settle.

On the other hand, V. Guba \cite{guba} constructed a finitely generated
almost cyclic group which is not cyclic. Guba refers in his paper to a
question of D. Anosov in the Kourovka Notebook which asks for the
construction of such an example.

Another application to topology arises as follows. If $x$ is a
conjugate generator for $G$, then the normal closure of $x$ is equal to
$G$. Kervaire \cite{kervaire} showed that $G$ is the fundamental group
of the complement of $S^n$ in $S^{n+2}$ for $n \ge 3$ if and only if the
abelianization of $G$ is $\Z$, $H_2(G) = 0$, and $G$ is the normal
closure of a single element.

\begin{question} If $G$ is the fundamental group of the complement of
$S^n$ in $S^{n+2}$ for $n \ge 3$, when is $G$ almost cyclic?
\end{question}

\section{Basic properties}

\begin{lemma} Let $G$ be almost cyclic, and let $Z(G)$ denote
the center of $G$. Then
$$Z(G) \subset \bigcap \left\{\langle x\rangle \mid x
   \quad\hbox{is a conjugate generator}\right\}.$$
\end{lemma}

\begin{proof} Let $g \in Z(G)$ and let $x$ be a conjugate generator. Then
$x^n = aga^{-1} = g$ for some $a \in G$ and some $n \in \Z$. This shows
that $Z(G) \subset \langle x\rangle$ for every conjugate generator
$x$.
\end{proof}

\begin{corollary} An abelian almost cyclic group is cyclic.
\end{corollary}

\begin{proof} If $G$ is abelian and almost cyclic and $x$ is a conjugate
generator for $G$, then $G = Z(G) \subset \langle x\rangle$. Hence, $G$
is cyclic, since subgroups of cyclic groups are cyclic.
\end{proof}

Likewise, a central subgroup of an almost cyclic group is cyclic.

\begin{proposition} The quotient of an almost cyclic group is
cyclic.
\end{proposition}

\begin{proof} Let $G$ be almost cyclic with conjugate generator $x$,
and let $N$ be a normal subgroup of $G$. I claim that $xN$ is a
conjugate generator for the quotient group $G/N$.

Let $gN \in G/N$. Find $a \in G$ such that $aga^{-1} = x^n$, where $n
\in \Z$. Then
$$(aN)(gN)(aN)^{-1} = (xN)^n.$$

Thus, $xN$ is a conjugate generator for $G/N$, and $G/N$ is almost
cyclic.
\end{proof}

The class of almost cyclic groups is not closed under products; for
example, $\Z_2 \times \Z_2$ is not almost cyclic. On the other hand, if
$G \times H$ is almost cyclic, the preceding result implies that $G$ and
$H$ are almost cyclic.

The following result is well-known (cf. \cite[Problem 1.9]{dixon}); a
proof is given for the convenience of the reader.

\begin{lemma} Let $G$ be a finite group, and let $H$ be a
proper subgroup of $G$. Then
$$\bigcup_{x\in G} xHx^{-1} \ne G.$$
\end{lemma}

\begin{proof} Let $n$ be the number of conjugates of $H$. If $n = 1$,
then $H$ is normal, and
$$\bigcup_{x\in G} xHx^{-1} = H \ne G.$$

Assume that $n \ge 2$. $G$ acts by conjugation on the set of conjugates
of $H$; the isotropy group of $H$ is the normalizer $N(H)$, so $n =
(G:N(H))$.

Since two conjugates have at least the identity in common, it follows
that
\begin{eqnarray*}
   \left|\bigcup_{x\in G} xHx^{-1}\right| & \le & nH - (n - 1) \\
   & = & (G:N(H))|H| - (n - 1) \\
   & \le & (G:H)|H| - (n - 1) \\
   & = & |G| - (n - 1) \\
   & < & |G| \\
\end{eqnarray*}
\end{proof}

\begin{proposition} A finite almost cyclic group is cyclic.
\end{proposition}

\begin{proof} Let $G$ be a finite almost cyclic group with conjugate
generator $x$. Suppose that $\langle x\rangle \ne G$. The preceding
lemma shows that
$$\bigcup_{g\in G} g\langle x\rangle g^{-1} \ne G.$$

Let $\displaystyle y \in G - \bigcup_{g\in G} g\langle x\rangle g^{-1}$.
Since $x$ is a conjugate generator, there is an element $a \in G$ and an
integer $n$ such that $a^{-1}ya = x^n$. Hence,
$$y = ax^na^{-1} \in \bigcup_{g\in G} g\langle x\rangle g^{-1}.$$

This contradiction implies that $\langle x\rangle = G$, and hence that
$G$ is cyclic.
\end{proof}

On the other hand, there are infinite almost cyclic groups which are not
cyclic. By a result of Higman, Neumann, and Neumann
\cite{higman-neumann-neumann}, a torsion-free group $G$ may be embedded
in a torsion-free group $G^*$ in which all nontrivial elements are
conjugate. In this case, any nontrivial element of $G^*$ is a conjugate
generator for $G^*$.

This result also shows that the class of almost cyclic groups is not
closed under passage to subgroups. For the Higman-Neumann-Neumann
construction applied to the rationals $\Q$ produces an almost cyclic
group $\Q^*$ with $\Q$ embedded as a subgroup. Since $\Q$ is abelian
but not cyclic, it is not almost cyclic.

\begin{lemma} Let $G$ be almost cyclic with conjugate generator $x$, and
let $H$ be a nontrivial subgroup of $G$. There is an element $g \in G$
such that $gHg^{-1} \cap \langle x\rangle \ne \{1\}$.
\label{normalgen}
\end{lemma}

\begin{proof} Let $h \in H$, $h \ne 1$. Find $g \in G$ such that
$ghg^{-1} = x^n$ for some $n \in \Z$. Note that $x^n \ne 1$, else $h =
1$. Hence, $ghg^{-1} = x^n$ is a nontrivial element of $gHg^{-1} \cap
\langle x\rangle$.
\end{proof}

\begin{lemma} Let $G$ be almost cyclic with conjugate generator $x$, and
let $H$ be a nontrivial normal subgroup of $G$. Let $n$ be the smallest
positive integer such that $x^n \in N$. Then
$$N \cap \langle x\rangle = \langle x^n\rangle.$$
\end{lemma}

\begin{proof} The preceding lemma implies that {\it some} conjugate
of $N$ intersects $\langle x\rangle$ nontrivially; since $N$ is normal,
it follows that $N$ intersects $\langle x\rangle$ nontrivially.

Let $n$ be the smallest positive integer such that $x^n \in N$. Then
$N \cap \langle x\rangle$ is cyclic (since it's a subgroup of $\langle
x\rangle$) with generator $x^n$.
\end{proof}

\begin{corollary} Let $G$ be almost cyclic with conjugate generator $x$
of prime order $p$. Then $G$ is simple.
\end{corollary}

\begin{proof} $\langle x\rangle$ is cyclic of prime order, so every
nontrivial element of $\langle x\rangle$ generates $\langle x\rangle$.

Let $N$ be a nontrivial normal subgroup of $G$. $N$ intersects $\langle
x\rangle$ nontrivially, so $\langle x\rangle \subset N$. Every element
of $G$ is conjugate to an element of $\langle x\rangle$, so every
element of $G$ is conjugate to an element of $N$. Since $N$ is normal,
every element of $G$ is contained in $N$, so $N = G$ and $G$ is
simple.
\end{proof}

In fact, in this case one can say more. The commutator subgroup $[G,G]$
is normal, so simplicity implies that $[G,G] = \{1\}$ or $[G,G] = G$. If
$[G,G] = \{1\}$, then $G$ is abelian, so it is cyclic --- in fact, $G
\approx \Z_p$. If $[G,G] = G$, $G$ is perfect --- in fact, an infinite
simple perfect group of exponent $p$.

\section{Solvable groups}

The main result of this section is that {\it solvable almost cyclic
groups are cyclic}. Consider the torsion case first.

Suppose that $G$ is a solvable almost cyclic group with conjugate
generator $x$ of finite order. Since $\langle x\rangle$ is finite, $G$
has only finitely many conjugacy classes. But a solvable group with
finitely many conjugacy classes is finite ([6], v.2, page 190). Hence,
$G$ is cyclic.

The torsion-free case is more work.

\begin{lemma} Let $G$ be almost cyclic with conjugate generator $x$. Let
$N$ be a normal subgroup of $G$ such that $G/N$ is cyclic. Then $G/N =
\langle x^kN\rangle$ for some $k \in \Z$.
\end{lemma}

\begin{proof} Let $yN$ be a generator of $G/N$. Then $gyg^{-1} = x^k$
for some $g \in G$ and $k \in \Z$, so $y = g^{-1}x^kg$ and $yN =
(g^{-1}x^kg)N$. Since $yN$ generates $G/N$, so does $(g^{-1}x^kg)N$; the
conjugate of a generator is also a generator, so $x^kN$ generates $G/N$
as well.
\end{proof}

\begin{lemma} Let $G$ be almost cyclic with conjugate generator $x$. Let
$N$ be a normal subgroup of $G$ such that $G/N$ is cyclic. Then:
\item{(a)} $N$ is almost cyclic.
\item{(b)} $G/N$ is {\it finite} cyclic.
\end{lemma}

\begin{proof} (a) By Lemma \ref{normalgen}, $N \cap \langle x\rangle =
\langle x^m\rangle$, where $x$ is the smallest positive power of $x$
contained in $N$. I claim that $x^m$ is a conjugate generator for $N$.

Let $y \in N$. Since $x$ is a conjugate generator for $G$,
$$gyg^{-1} = x^p \quad\hbox{for some}\quad g \in G, \quad
   p \in \Z.$$

Since $N$ is normal, $gyg^{-1} = x^p \in N \cap \langle x\rangle =
\langle x^m\rangle$, so $gyg^{-1} = (x^m)^q$ for some $q \in \Z$. By the
preceding lemma, $G/N = \langle x^kN\rangle$ for some $k \in \Z$. Thus,
$g = (x^k)^rz$ for some $r \in \Z$ and $z \in N$. Then
$$x^{kr}zyz^{-1}x^{-kr} = (x^m)^q, \quad\hbox{and hence}\quad
   zyz^{-1} = (x^m)^q.$$

That is, $y$ is conjugate by an element $z \in N$ to a power of $x^m$.
This proves the claim, and the lemma.

(b) Since $x^m \in N$, and since $G/N = \langle x^kN\rangle$, it's clear
than $m$ is an exponent for $G/N$. Therefore, $G/N$ is finite
cyclic.
\end{proof}

\begin{theorem} Let $G$ be a solvable torsion-free almost cyclic
group. Then $G$ is cyclic.
\end{theorem}

\begin{proof} Induct on the length of a solvable series for $G$:
$$G = G_0 \triangleright G_1 \triangleright \cdots \triangleright
   G_n = \{1\}.$$

If $n = 1$, then $G$ is abelian, so it is cyclic.

Suppose that $G$ has a solvable series of length $n > 1$, and suppose
that the result is true for solvable torsion-free almost cyclic groups
with solvable length less than $n$. Let
$$G = G_0 \triangleright G_1 \triangleright \cdots \triangleright
   G_n = \{1\}$$

be a solvable series for $G$, and consider the exact sequence
$$1 \longrightarrow G_1 \longrightarrow G \longrightarrow
   \dfrac{G}{G_1} \longrightarrow 0.$$

By assumption, $G/G_1$ is abelian; since it's a quotient of an almost
cyclic group, it's almost cyclic. Therefore, $G/G_1$ is cyclic --- in
fact, {\it finite} cyclic, by the preceding lemma..

The preceding lemma also implies that $G_1$ is almost cyclic. Thus,
$G_1$ is a solvable torsion-free almost cyclic group with a solvable
series of length $n - 1$. By induction, $G_1$ is cyclic; hence, it is
infinite cyclic, since it is torsion-free.

Finally, $G$ is torsion-free, and it has an infinite cyclic subgroup of
finite index, so it is infinite cyclic (\cite[page 96]{brown}).
\end{proof}

\section{Almost cyclic one-relator groups}

As noted above, in \cite{guba}, Guba constructed a finitely generated
almost cyclic group that is not cyclic.

\begin{question} Is every finitely presented almost cyclic group cyclic?
\end{question}

This problem is likely very difficult. The aim of this section is to
give an affirmative answer for the special case of a one-relator group.

\begin{proposition} Let $G = \langle X\mid r\rangle$ be an almost cyclic
one-relator group. If $r$ is a proper power, then $G$ is finite cyclic.
\end{proposition}

\begin{proof} Let $x$ be a conjugate generator for $G$. Since $r$ is
a proper power, it can be written in the form $r = s^m$, where $s$ is a
word in $X \cup X^{-1}$ and $m > 1$ is maximal. Let $\overline{s}$ be
the image of $s$ in $G$. Since $x$ is a conjugate generator,
$$g\overline{s}g^{-1} = x^n \quad\hbox{for some}\quad
   g \in G, \quad n \in \Z.$$

Then $g\overline{s}^mg^{-1} = x^{mn}$, but $\overline{s}^m = 1$ in $G$.
Hence, $x^{mn} = 1$, and $G$ is a torsion group.

However, Fischer, Karrass, and Solitar \cite{fischer-karrass-solitar}
have shown that every one-relator group is virtually torsion-free. Since
$G$ is a torsion group, $G$ must be finite. Therefore, $G$ is finite
cyclic.
\end{proof}

A result of Lyndon \cite{lyndon} implies that if the relator $r$ in a
one-relator group $G = \langle X\mid r\rangle$ is not a proper power,
then $G$ is torsion-free.

The next result is not strictly necessary, but it will simplify the
notation in the theorem that follows.

\begin{lemma} Let $G = \langle X \mid r\rangle$ be a one-relator almost
cyclic group. Then $|X| \le 2$.
\end{lemma}

\begin{proof} The abelianization $\dfrac{G}{[G,G]}$ is the quotient
of an almost cyclic group, so it's almost cyclic. Since it's abelian,
it's cyclic.

On the other hand, if $|X| > 2$, then $\dfrac{G}{[G,G]}$ is a direct
product of cyclic groups with at least two infinite cyclic
factors. Therefore, $|X| \le 2$.
\end{proof}

The same argument gives an easy proof that the only almost cyclic free
group is $\Z$.

\begin{theorem} Let $G = \langle X \mid r\rangle$ be a one-relator
group, and suppose $G$ is almost cyclic with conjugate generator
$x$. Then $G$ is cyclic.
\end{theorem}

\begin{proof} The result is obvious if $G$ is generated by a single
element. By the preceding lemma, the only case left to consider is $|X|
= 2$. Suppose then that $G = \langle t, u \mid r\rangle$ with $r \in
\langle t,u\rangle$. As usual, let $\sigma_s(w)$ denote the exponent sum
of the generator $s$ in the word $w$. If $w$ is a word in the free group
$\langle t,u\rangle$, then $\overline(w)$ will denote the image in $G$,
regarded as the quotient $\dfrac{\langle t,u\rangle}{\langle
r\rangle_n}$.

{\bf Case 1.} $r$ has zero exponent sum on one of the generators.

Without loss of generality, suppose that $\sigma_t(r) = 0$. Since $x$ is
a conjugate generator for $G$,
$$a\overline{t}a^{-1} = x^m \quad\hbox{and}\quad
   b\overline{u}b^{-1} = x^n$$

for some $a, b \in G$ and $m, n \in \Z$. Then
$$a\overline{t}^na^{-1} = x^{mn} = b\overline{u}^mb^{-1}.$$

Let $a_0$, $b_0$ be lifts of $a$, $b$, respectively, to $\langle
t,u\rangle$. Then
$$a_0t^na_0^{-1} = b_0u^mb_0^{-1} \mod{\langle r\rangle_n}
   \quad\hbox{in}\quad \langle t,u\rangle.$$

An element of $\langle r\rangle_n$ is a product of conjugates of $r$ or
$r^{-1}$. Since $\sigma_t(r) = 0$, the same is true of every element of 
$\langle r\rangle_n$. Therefore,
$$\sigma_t(a_0t^na_0^{-1}) = \sigma_t(b_0u^mb_0^{-1}),
   \quad\hbox{and}\quad n = 0.$$

Thus, $b\overline{u}b^{-1} = 1$, and $\overline{u} = 1$. Therefore, $G$
is cyclic.

{\bf Case 2.} $\sigma_t(r) = p \ne 0$ and $\sigma_u(r) = q \ne 0$.

Every one-relator group with at least two generators admits a
one-relator presentation in which the relator has zero exponent sum on
at least one relator \cite[Lemma 11.9]{lyndon-schupp}. Therefore, this
case reduces to the previous one.
\end{proof}

\begin{question} Is every almost cyclic group satisfying a small
cancellation condition cyclic?
\end{question}

\bibliography{almcyc}
\bibliographystyle{plain}

\end{document}